\documentclass{amsart}[12pt]
\usepackage{graphics}
\usepackage{verbatim}

\newtheorem{thm}{Theorem}
\newtheorem{cor}{Corollary}
\newtheorem{prop}{Proposition}
\newtheorem{lemma}{Lemma}
\theoremstyle{definition} \newtheorem{define}{Definition}
\theoremstyle{remark} 
\theoremstyle{remark}  \newtheorem*{example}{Example}
\theoremstyle{remark}  \newtheorem{exercise}{Exercise}

\newcommand{\xhelp}[1]{\textbf{$\diamondsuit\diamondsuit$ Fix me: #1}}

\newcommand{\RR}{\mathbb{R}}
\newcommand{\CC}{\mathbb{C}}
\newcommand{\ZZ}{\mathbb{Z}}
\newcommand{\ff}{\vec{f}}
\newcommand{\vv}{\vec{v}}
\newcommand{\ww}{\vec{w}}
\newcommand{\id}{\mathbf{1}}
\newcommand{\Meas}{\operatorname{\textrm{Meas}}}
\renewcommand{\ker}{\operatorname{\textrm{Ker}}}
\newcommand{\im}{\operatorname{\textrm{Im}}}
\newcommand{\vect}{\operatorname{\textrm{Vec}}}
\newcommand{\inter}[1]{{#1}^{\circ}}
\newcommand{\mucan}{{\mu_{\rm can}}}
\newcommand{\aff}{\operatorname{A(\Gamma)}}
\newcommand{\smooth}{\operatorname{S(\Gamma)}}
\newcommand{\smoothmu}{\operatorname{S_{\mu}(\Gamma)}}

\title[Metrized graphs]
{Metrized graphs, electrical networks, and Fourier analysis}

\author{Matt Baker}
\address{Matt Baker\\
School of Mathematics\\
Georgia Institute of Technology\\
Atlanta, GA 30332-0160\\
USA}
\email{mbaker@math.gatech.edu}

\author{Xander Faber}
\address{Xander Faber\\
Department of Mathematics\\
Columbia University\\
New York, NY  10027\\
USA}
\email{xander@math.columbia.edu}

\thanks{The authors would like to thank Robert Rumely for numerous
  useful conversations about metrized graphs.  They also thank their
  wives for tolerating being ignored while this article was being
  written.  Finally, they would like to thank the participants of the
  Summer 2003 REU program ``Analysis on Metrized Graphs'' for their
  dedication and enthusiasm: Maxim Arap, Jake Boggan, Rommel Cortez,
  Crystal Gordon, Kevin Mills, Kinsey Rowe, and Phil Zeyliger.  The
  first author's work was supported by NSF grant DMS 0070736. The
  second author's work was supported by the VIGRE grant at Columbia
  University (DMS 9810750).}

\begin{document}

\begin{abstract}
A metrized graph is a finite weighted graph whose edges are thought of as line segments.
In this expository paper, we study the Laplacian operator on a metrized graph and some
important functions related to it, including the ``j-function'', the effective
resistance, and eigenfunctions of the Laplacian. We discuss the relationship between
metrized graphs and electrical networks, which provides some physical intuition for the
concepts being dealt with. We also discuss the relation between the Laplacian on a
metrized graph and the combinatorial Laplacian matrix. We introduce the``canonical
measure'' on a metrized graph, which arises naturally when considering the Laplacian of
the effective resistance function. Finally, we discuss a generalization of classical
Fourier analysis which utilizes eigenfunctions of the Laplacian on a metrized graph.
During the course of the paper, we obtain a proof of Foster's network theorem and of the
identity $\min\{x,y\} = 8 \sum_{n \geq 1 \textrm{ odd}} \sin(n\pi x/2) \sin(n\pi y/2)/
\pi^2 n^2$, for $0 \leq x,y \leq 1$.
\end{abstract}

\maketitle

\section{An informal discussion}\label{discussion}

Graphs are usually considered to be discrete objects, so issues of continuity and
differentiability don't typically appear in graph theory texts.  Here is a picture of a
graph:

\begin{figure}[!h]
  \scalebox{.6}{\includegraphics{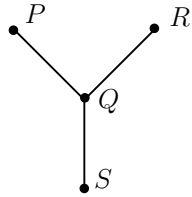}}

  \caption{A graph with four vertices and three edges.}

  \label{Flux Capacitor}
\end{figure}

The picture is misleading, in the sense that we might want to believe that the edge $PQ$
is a line segment comprising a continuum of points. However, edges in graph theory are
merely formal connections between vertices; they don't ``contain points.''  But is there
something {\em wrong} with thinking of $PQ$ as a line segment?  
The notion of a \textit{metrized graph} gives meaning to points between the vertices
while retaining the salient combinatorial features of the graph.  In a broader sense,
metrized graphs will unite the discrete and the continuous for us.

The basic idea of a metrized graph is simple: identify each edge of a finite
graph with a line segment and define the distance between two points of the
graph to be the length of the shortest path connecting them. We will provide more
details on this definition in \S\ref{MetrizedEquivalence}.

Metrized graphs appear in the literature of several areas of science and mathematics. For
example: in number theory, they are used to study arithmetic intersection theory on
algebraic curves (see \cite{CR}, \cite{Zhang}); in mathematical biology, they are used to
study neuron transmission (see \cite{Ni}); they are also used in physics, chemistry, and
engineering (under the names \textit{metric graphs}, \textit{quantum graphs}, and
$c^2$\textit{-networks}) as wave-propagation models (see \cite{Ku}).

In \S\ref{Continuous Laplacian} we define a Laplacian operator on a metrized graph which
is closely related to the Laplacian matrix (or Kirchhoff matrix) associated to a finite
graph---see \S\ref{Laplacian Matrix} for precise statements about this connection.  Using
a theory of eigenfunction expansions on a metrized graph that generalizes classical
Fourier analysis on the circle, we will be able to prove intriguing series identities
such as
\[
\min\{x,y\} = 8 \sum_{n \geq 1 \textrm{ odd}} \frac{\sin \left(\frac{n \pi x}{2}\right)
\sin \left(\frac{n \pi y}{2}\right)}{ \pi^2 n^2}, \qquad 0 \leq x,y \leq 1.
\]
See \S\ref{Fourier Analysis} for more details.\footnote{Metrized graphs can be viewed as
one-dimensional Riemannian manifolds with singularities, and from this point of view the
Laplacian on a metrized graph is a nontrivial but computationally accessible variant of
the Laplacian on a higher-dimensional Riemannian manifold.}

There is a well-known and useful interplay between the theories of finite graphs and
resistive electrical networks (see e.g., \cite[Ch.~II,IX]{Bo}). This relationship extends
beautifully to the setting of metrized graphs (cf. \S\ref{Continuous
Laplacian},\ref{jFunctionSection}). For example, a theorem of Foster from 1949 (see
\cite{Fo}) asserts that
\[
\sum_{\textrm{edges } e} \frac{r(e)}{L_e} = \#V - 1,
\]
where $r(e)$ is the effective resistance in the electrical network between the endpoints
of the edge $e$, $L_e$ is the resistance along the edge $e$, and $\#V$ is the number of
nodes (vertices) in the network. In \S\ref{Canonical Measure Section} we give a proof of
Foster's theorem using the ``canonical measure'' on a metrized graph.

The theory of electrical networks is itself closely related to the theory of random walks
on graphs.  We will not touch upon the connection with random walks in this paper, but we
refer the interested reader to the delightful monograph \cite{DS}. There is a nice proof
of Foster's theorem using random walks in \cite{Te} (see also \cite[Theorem 25, Exercise
23, Chapter IX]{Bo}).

This article is a follow-up to the 2003 summer REU on metrized graphs held at the
University of Georgia and run by the first author and Robert Rumely. The participating
students' enthusiasm for the subject convinced us that a broader audience might
appreciate an introduction to the ideas involved. Further information about the REU, its
organizers and participants, and the research they performed can be found at
\verb+http://www.math.uga.edu/~mbaker/REU/REU.html+

In keeping with the spirit of discovery that spawned this article, we have included a
number of exercises to clarify the text or extend the ideas presented. We have also
strived to keep the exposition as self-contained as possible with the hopes that it will
inspire further students toward this subject.


\section{Metrized graphs versus weighted graphs} \label{MetrizedEquivalence}

There is a bijective correspondence between metrized graphs and equivalence
classes of finite weighted graphs. In this section we give an overview of
this correspondence, leaving many of the details to the reader.
See \cite{BR} for more detailed proofs of the assertions made in this section.

\begin{define}
For the purposes of this paper, we define a \textit{weighted graph} $G$ to
be a finite, connected graph with vertex set $V(G) = \{ v_1,\ldots,v_n \}$,
edge set $E(G) = \{e_1,\ldots,e_m\}$, and a collection of positive weights
$\{ w_{e_1} \ldots, w_{e_m} \}$ associated to the edges of $G$. Further, we
require that $G$ have no loop edges or multiple edges. The \textit{length}
of the edge $e$ is defined to be $L_e=1/w_e$.
\end{define}

\begin{figure}[!ht]

    \begin{picture}(200,65)(0,0)
        \put(-15,0){\scalebox{0.6}{\includegraphics{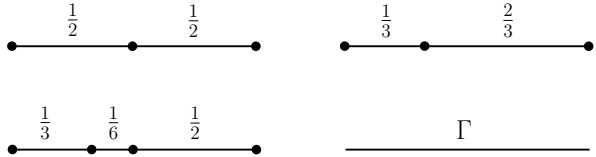}}}

    \end{picture}

\caption{Each of the three weighted graphs displayed is a model of the metrized graph
$\Gamma$, a segment of length $1$. The weighted graphs are all distinct, but they lie in
the same equivalence class. The lower left weighted graph is a common refinement of the
upper left and upper right.} \label{Equivalent Graphs}
\end{figure}

In classical graph theory, one associates weights to the edges of a graph. When studying
metrized graphs, it makes more sense to work with lengths, since distance is the
fundamental notion in a metric space. We will henceforth indicate {\em lengths} in our
figures (e.g., Figure~\ref{Equivalent Graphs}).

A weighted graph $G$ gives rise to a metric space $\Gamma$ in the following
way. To each edge $e$, associate a line segment of length $L_e$, and identify
the ends of distinct line segments if they correspond to the same vertex of
$G$. The points of these line segments are the points of $\Gamma$. We call
$G$ a $\textit{model}$ for $\Gamma$.
The distance between two points $x$ and $y$ in
$\Gamma$ is defined to be the length of the shortest path between them, where
the length of a path is measured in the usual way along the line segments
traversed.  (A path between distinct points always exists because $G$ is
connected.)

\begin{exercise}
Show that this notion of distance
defines a metric on $\Gamma$ (which we call the \textit{path metric}).
\end{exercise}

\begin{figure}[!ht]
  \begin{picture}(200,55)(0,0)
    \put(-30,10){\scalebox{0.55}{\includegraphics{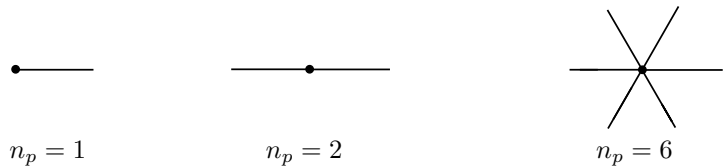}}}
    \put(-30,0){$n_p=1$}
    \put(67,0){$n_p=2$}
    \put(192,0){$n_p=6$}

  \end{picture}

  \caption{Three examples of star-shaped sets and their valences.}
  \label{Star-Shaped Sets}
\end{figure}

The space $\Gamma$,
endowed with the path metric, is called a \textit{metrized graph}.
Here's a more abstract definition, taken from \cite{Zhang}:

\begin{define} \label{AbstractMetrizedGraph}
A \textit{metrized graph} $\Gamma$ is a compact, connected metric space such
that each $p \in \Gamma$ has a
neighborhood $U_p$ isometric to a star-shaped set of valence $n_p \geq 1$, endowed with
the path metric (see Figure~\ref{Star-Shaped Sets}). To be precise, a star-shaped set of
valence $n_p$ is a set of the form
\[ S(n_p,r_p) \ = \ \{ z
\in \CC : z = t e^{k \cdot 2 \pi i/n_p} \ \text{for some $0 \le t < r_p$ and some $k \in
\ZZ$} \}. \]
\end{define}

\begin{exercise}
Check that the metric space $\Gamma$ arising from a
weighted graph $G$ satisfies the abstract definition
(Definition \ref{AbstractMetrizedGraph}) of a metrized graph.
\end{exercise}

The points $p \in \Gamma$ with valence different from 2 are precisely those where $\Gamma$
fails to look locally like an open interval, and the compactness of $\Gamma$ ensures
that there are only finitely many such points. Let $V(\Gamma)$ be any
finite, nonempty subset of $\Gamma$ such that:
\begin{itemize}
\item $V(\Gamma)$ contains all of the points with $n_p \not= 2$.
(This implies that $\Gamma \setminus V(\Gamma)$ is a finite, disjoint union of subspaces $U_i$
isometric to open intervals.)
\item For each $i$, the topological closure $\overline{U}_i$ of $U_i$ in $\Gamma$ is isometric
 to a line segment (as opposed to a circle).  We call $e_i =
 \overline{U}_i$ a \textit{segment} of $\Gamma$.
\item For each $i \neq j$,
 $e_i \cap e_j = \emptyset$ or $\{p \}$, where $p$
 is an endpoint of both $e_i$ and $e_j$.

\end{itemize}

Any finite set $V(\Gamma)$ satisfying these conditions will be called a
\textit{vertex set} for $\Gamma$, and the elements of $V(\Gamma)$ will
be called \textit{vertices} of $\Gamma$.

\begin{exercise}

\begin{itemize} \item[] \hspace{1in}
\item[(a)] Prove that a vertex set for $\Gamma$ always exists.
\item[(b)] Let $\Gamma$ be a circle. Show that any set consisting of only one or two points of
$\Gamma$ cannot be a vertex set.
\end{itemize}
\end{exercise}

It should be remarked that $V(\Gamma)$ is not unique.  For example, if $\Gamma$ is a
circle, then any choice of three distinct points of $\Gamma$ is a vertex set. The choice
of a vertex set $V(\Gamma)$ determines a finite set $\{ e_i \}$ of segments of $\Gamma$.
The endpoints of each segment $e_i$ are vertices of $\Gamma$. We emphasize that the
segments of $\Gamma$ depend on our choice of a vertex set.

Given a metrized graph $\Gamma$, our next task will be to find a weighted graph $G$ that
serves as a model for $\Gamma$ as above. Pick a vertex set $V(\Gamma)$ for $\Gamma$.
Define a graph $G$ with vertices indexed by $V(\Gamma)$, and join two distinct vertices
$p$ and $q$ of $G$ by an edge if and only if there exists a segment of $\Gamma$ with
endpoints $p$ and $q$. (So edges of $G$ correspond to segments of $\Gamma$.) Define the
length of the edge joining $p$ to $q$ to be the length of the segment $e$.  Then $G$ is a
weighted graph, with weights given by the reciprocals of the lengths; our definition of
$V(\Gamma)$ guarantees that $G$ has no multiple edges or loop edges.  Moreover, if we
construct the metrized graph associated to $G$, it is easily seen to be isometric to
$\Gamma$.

Different choices of a vertex set $V(\Gamma)$ yield distinct weighted graphs in the above
construction. Write $G \sim G'$ if the two weighted graphs $G, G'$ admit a \textit{common
refinement}, where we \textit{refine} a weighted graph by subdividing its edges in a
manner that preserves total length (see Figure~\ref{Equivalent Graphs}). This provides an
equivalence relation on the collection of weighted graphs, and one can check that two
weighted graphs are equivalent if and only if they give rise to isometric metrized
graphs.

Having established this correspondence, we are now free to fix a particular model of a
metrized graph, without worrying that we've lost some degree of generality in doing so.
This will be especially convenient in the next section, when we work with functions on a
metrized graph that are ``nice'' outside of some vertex set.


\section{The Laplacian on a metrized graph} \label{Continuous Laplacian}

Our goal in this section is to motivate and define the Laplacian of a function on a
metrized graph.  The Laplacian on a metrized graph is a hybrid between the Laplacian on
the real line (i.e., the negative of the second derivative) and the discrete Laplacian
matrix studied in graph theory (cf. \S\ref{Laplacian Matrix}).

Choose a vertex set $V(\Gamma)$ for the metrized graph $\Gamma$. Let $p$ be
a non-vertex point of $\Gamma$, and suppose $e$ is a segment of length $L$
containing $p$. Parametrize $e$ by an isometry $s_e : [0,L]\to e$ so that we
have a real coordinate $t \in [0,L]$ to use for describing points of the
segment. We say that $f$ is differentiable at $p$ if the quantity
$\frac{d}{dt}f(s_e(t))|_{s_e(t)=p}$ exists. There is precisely one other
parametrization of this sort, namely $u_e(t)=s_e(L-t)$. The chain rule shows
that
\[
 \frac{d}{dt}f(u_e(t))\Big|_{u_e(t)=p}=
-\frac{d}{dt}f(s_e(t))\Big|_{s_e(t)=p}.
\]
Hence the value of the derivative of $f$ at $p$ depends on the
parametrization, but only up to a sign. Picking one of the two
parametrizations for a segment can be thought of as choosing an
\textit{orientation} for the segment, and we will use the two concepts
interchangeably.

We can similarly determine if $f$ is $n$ times differentiable at $p$ by looking at the existence
of the quantity $\frac{d^n}{dt^n} f(s_e(t))|_{s_e(t)=p}$.

\begin{exercise}
Show that the second derivative $f''(p)$, when it exists, is
well-defined independent of the choice of an orientation for $e$.
\end{exercise}

We also require a notion of differentiability that makes sense at the vertices.
The abstract definition of a metrized graph tells us that each point $p \in
\Gamma$ has a neighborhood isometric to a star-shaped set with $n_p \geq 1$
arms. Thus there are $n_p$ directions by which a path in $\Gamma$ can
leave $p$.  To each such direction, we associate a formal unit vector $\vv$, and we
write $\vect(p)$ for the collection of all $n_p$ directions at $p$. We make
this convention so that we can write $p+\varepsilon\vv$ for the point of
$\Gamma$ at distance $\varepsilon$ from $p$ in the direction $\vv$ for
sufficiently small $\varepsilon > 0$.

\begin{define} \label{directderiv}
Given a function $f:\Gamma \to \RR$, a point $p \in \Gamma$, and a direction $\vv \in
\vect(p)$, the \textit{derivative of $f$ at $p$ in the direction $\vv$}, written
$D_{\vv}f(p)$, is given by
\[
D_{\vv}f(p) = \lim_{\varepsilon\to 0^+} \frac{f(p+\varepsilon\vv) -f(p)}{\varepsilon},
\]
provided this limit exists. This will also be called a \textit{directional derivative}.
\end{define}

\begin{exercise}
\label{SigmaFiniteExercise}
Given a function $f: \Gamma \to \RR$ and a point
$p \not\in V(\Gamma)$ at which $f$ is differentiable, show that the two directional derivatives of
$f$ at $p$ exist and sum to zero.
[{\bf Hint:} Parametrize the segment containing $p$ and do
the calculation explicitly.]
\end{exercise}

Here is the class of functions on which we intend to apply our Laplacian:
\begin{define} \label{smooth}
Define $\smooth$ to be the class of all continuous functions $f: \Gamma \to \RR$ for
which there exists a vertex set $V_f(\Gamma)$ (with corresponding segments $e_i$) such
that
\begin{itemize}
\item[(i)] $D_{\vv}f(p)$ exists for each $p \in \Gamma$ and each $\vv \in
\vect(p)$,
\item[(ii)] $f$ is twice continuously differentiable on the
  interior of each segment $e_i$, and
\item[(iii)] $f''$ is bounded on the interior of each segment $e_i$.
\end{itemize}
We call $\smooth$ the class of \textit{piecewise smooth functions} on $\Gamma$. (This is,
of course, a small abuse of terminology as these functions need not be infinitely
differentiable away from the vertices.)
\end{define}

\begin{exercise}
Show that hypotheses (ii) and (iii) imply hypothesis (i), and that
hypothesis (i) already implies that $f$ is continuous.
\end{exercise}

We now define the Laplacian operator on a metrized graph. A conceptual obstacle to
overcome is that the Laplacian of a function $f \in \smooth$ is a \textit{bounded, signed
measure}\footnote{We could also work with complex-valued functions, and then the
Laplacian would be a complex measure.} on $\Gamma$, not a function. For readers
unfamiliar with the notion of measure, we give a brief working definition in just a
moment.

\begin{define} \label{Laplacian}
The \textit{Laplacian} of a function $f \in \smooth$ is given by the
measure
\[
\Delta f = -f''(x) dx - \sum_{p \in \Gamma} \sigma_p(f) \delta_p,
\]
where $\sigma_p(f) = \sum_{\vv \in \vect(p)} D_{\vv}f(p)$, $dx$ denotes the Lebesgue
measure on $\Gamma$, and $\delta_p$ is the Dirac measure (unit point
mass) at $p$.
\end{define}
By Exercise~\ref{SigmaFiniteExercise}, the sum $\sum_p \sigma_p(f)$
is actually finite as $\sigma_p(f) = 0$ for any $p$ not in $V_f(\Gamma)$. Also, $f''$ is
well-defined away from the vertices in $V_f(\Gamma)$, so $\Delta f$ is independent of
segment orientations. To perform computations with $\Delta f$, however, we will need to
choose parametrizations.

Let's define some notation to make what lies ahead a little easier.  Choose a model for
$\Gamma$, parametrize each segment $e$ of $\Gamma$ by $s_e:[0,L_e] \to e$, and for $f:
\Gamma \to \RR$ define $f_e:[0,L_e] \to \RR$ by $f_e = f \circ s_e$.

Now for our working definition of measure.  Intuitively, a measure is
something we can integrate functions against. For our purposes, then, a measure on
$\Gamma$ will be an expression of the form
\[
\mu = \sum_{\textrm{segments } e} g_e(t)dt|_e + \sum_{i=1}^n c_i \delta_{p_i},
\]
where $g_e : (0,L_e) \to \RR$ is continuous and bounded, $c_i \in \RR$, and
$p_1,\ldots,p_n$ are points of $\Gamma$. To integrate a continuous function $f: \Gamma
\to \RR$ against the measure $\mu$, define
\[
\int_{\Gamma} f(x)d\mu(x) = \sum_{\textrm{segments } e} \left\{ \int_{0}^{L_e}
f_e(t)g_e(t) \, dt\right\} + \sum_i c_i f(p_i).
\]

A measure of the form $\sum g_e(t) dt|_e$ is called a
\textit{continuous measure}, and a measure of the form
$\sum_{i=1}^n c_i \delta_{p_i}$ is called a \textit{discrete
measure}.

If $g : \Gamma \to \RR$ is a function such that $g \circ s_e(t) = g_e(t)$ for all
segments $e$ of $\Gamma$ and all $t \in (0,L_e)$, we will usually write $g(x)dx$ instead
of $\sum_e g_e(t) dt|_e$.

The {\em total mass} of a measure $\mu$ is defined to be $\int_\Gamma \id(x) d\mu(x)$,
where $\id$ denotes the constant function with value 1.

\begin{example} \label{Laplacian Calculation}
Consider the metrized graph $\Gamma$ modelled in Figure~\ref{Flux Laplacian Example}.
Define a function on $\Gamma$ by
\[
f_e(t) = \begin{cases} t+1, & e = PQ \\ 3(t+\frac{1}{2}), & e= QS \\
t^2 + \frac{1}{2}, & e = RQ.\end{cases}
\]
Then $\Delta f = -2dx|_{RQ} - \delta_P + 3 \delta_S$.  Note that
$\Delta f$ has total mass zero; we will see shortly that this is not
an accident.

\begin{figure}[!ht]
  \begin{picture}(180,75)(0,0)
    \put(-20,0){\scalebox{0.6}{\includegraphics{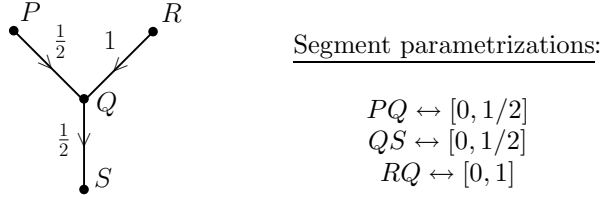}}}

    \put(80,30){
        \begin{tabular*}{126pt}{c}
        \underline{Segment parametrizations}: \\
        \\
        $PQ \leftrightarrow [0,1/2]$ \\
        $QS \leftrightarrow [0,1/2]$ \\
        $RQ \leftrightarrow [0,1]$
        \end{tabular*}
    }

  \end{picture}

  \caption{A model of a metrized graph and its segment parametrizations. The arrows in the diagram
    indicate the directions in which the segment parametrizations increase.}
  \label{Flux Laplacian Example}
\end{figure}
\end{example}

Our next goal is to put together some important facts about the Laplacian.
It will be more convenient to work with an alternate formulation,
though. For simplicity, we write $f_e'(0)$ for the right-hand derivative of
$f_e$ at $0$ (as the limit only makes sense from one side). Similarly, write
$f_e'(L_e)$ for the left-hand derivative at $L_e$. If $p$ and $q$ are the
endpoints of the segment $e$, with $s_e(0) = p$ and $s_e(L_e)=q$, we'll say that
$e$ \textit{begins} at $p$ and \textit{ends} at $q$. If $\vv \in \vect(p)$ and $\ww \in
\vect(q)$ are the directions pointing inward along $e$, then
\[
D_{\vv}f(p) = f_e'(0), \qquad D_{\ww}f(q) = -f_e'(L_e).
\]

Observe that $-\sigma_p(f)$ counts $-f_e'(0)$ for each segment $e$ beginning at $p$, and
it counts $f_e'(L_e)$ for each segment $e$ ending at $p$. Thus $\Delta
f$ can be written
\begin{equation*}\tag{$\ast$}
\label{AlternateDefinition}
\Delta f =
\sum_{\textrm{segments } e}
\left\{
-f_e''(x) \, dx|_e +
f_e'(L_e)\delta_{s_e(L_e)} -
f_e'(0)\delta_{s_e(0)} \right\}.
\end{equation*}
In particular, we now see that the contribution of the segment $e$ to the Laplacian is
$-f_e''(x)dx|_e +f_e'(L_e)\delta_{s_e(L_e)} - f_e'(0)\delta_{s_e(0)}$. This measure is
independent of the choice of parametrization of $e$, but
it is necessary to choose a parametrization to write it down.

\begin{thm}[Self-adjointness of $\Delta$] \label{Laplacian Fact 1}
Suppose $f, g \in \smooth$. Then
\[
\int_{\Gamma} f \Delta g = \int_{\Gamma} g \Delta f =
    \int_{\Gamma} f'(x)  g'(x) \, dx.
\]
\end{thm}

\begin{proof}
Choose a model for $\Gamma$ with vertex set $V(\Gamma) = V_f(\Gamma) \cup V_g(\Gamma)$.
Then $f''$ is continuous on the interior of each segment of $\Gamma$, and the directional derivatives of
$f$ and $g$ exist for all vertices in $V(\Gamma)$ (see Definition~\ref{smooth}). Choose
parametrizations for each segment of $\Gamma$ and define $f_e$ and $g_e$ as
before. Using integration by parts, we obtain:
\begin{align*}
\int_{\Gamma} g \Delta f
&= \sum_{e}
\left\{ f_e'(L_e)
g_e(L_e) - f_e'(0) g_e(0) - \int_{0}^{L_e} g_e(t) f_e''(t) \, dt \right\}\\
&= \sum_{e}
\int_{0}^{L_e} f_e'(t)g_e'(t) \, dx = \int_{\Gamma} f'(x)g'(x) \, dx.
\end{align*}

The rest of the result follows by symmetry.
\end{proof}


\begin{cor} \label{Laplacian Fact 2}
If $f \in \smooth$, then $\Delta f$ has total mass~$0$.
\end{cor}

\begin{proof}
Set $g = \id$ in the statement of Theorem~\ref{Laplacian Fact 1}. Then
\[\int_{\Gamma}
\id \cdot \Delta f = \int_{\Gamma} f'(x) \id'(x) \, dx = \int_\Gamma 0 \, dx = 0.\]
\end{proof}

Before stating the next result about the Laplacian, we need to define another useful
class of functions:
\begin{define} \label{affine}
Define $\aff$ to be the subclass of functions $f \in \smooth$ such that for
each oriented segment $e$ of $\Gamma$, there exist real constants $A_e,B_e$ so
that $f_e(t) = A_e t+B_e$ for $t \in [0,L_e]$. A function in $\aff$ is
called \textit{piecewise affine}.
\end{define}

\begin{exercise}
\label{AffineDetermined}
Show that a function $f \in \aff$ is completely
determined by its values on a vertex set $V_f(\Gamma)$ for $f$.
[{\bf Hint:} Use linear interpolation.]
\end{exercise}

\begin{exercise}
\label{AffineDiscrete}
Show that if $f \in \smooth$, then $f$ is piecewise affine
if and only if $\Delta f$ is a discrete measure.
\end{exercise}

The next result is a graph-theoretic analogue of the second derivative test from
calculus.

\begin{thm}[The Maximum Principle] \label{Max Principle}
Suppose $f \in \aff$ is nonconstant. Then $f$ achieves its maximum
value on $\Gamma$ at a vertex $p \in V_f(\Gamma)$ for which
$\sigma_p(f) < 0$.
\end{thm}

\begin{proof}
It is easy to see that the function $f$ must take on its maximum
value at a vertex $p \in V_f(\Gamma)$.
Moreover, since $f$ is nonconstant, we may select $p$ so that $f$ decreases
along some segment $e_0$ having $p$ as an endpoint.  (This uses the
fact that $\Gamma$ is connected.)  Re-parametrize each
segment $e$ having $p$ as an endpoint, if necessary, so that $s_e(0)=p$, where $s_e:[0,L_e] \to
e$. Then $\sigma_p(f) = \sum f_e'(0)$, where the summation is over all segments $e$
beginning at $p$. Each of the slopes $f_e'(0)$ must be non-positive;
otherwise $f$ would grow along $e$, violating the fact that $f$ is maximized
at $p$. We know $f_{e_0}'(0) < 0$ since $f$ decreases along $e_0$. Hence
$\sigma_p(f) < 0$, which completes the proof.
\end{proof}

\begin{thm} \label{Laplacian Fact 3}
Suppose $f, g \in \smooth$. If $\Delta f = \Delta g$ and $f(p)=g(p)$ for some $p \in
\Gamma$, then $f \equiv g$.
\end{thm}

\begin{proof}
If $h=f-g$, then $\Delta h = 0$.  By Exercise~\ref{AffineDiscrete}, $h \in \aff$.
As $\Delta h = 0$, it follows from the Maximum Principle that $h$ is constant.
The hypothesis that $f(p)=g(p)$ for some $p$
now implies that $h \equiv 0$, so that $f \equiv g$ as desired.
\end{proof}
Note in particular that if $f \in \smooth$ is {\em harmonic} (i.e., $\Delta f = 0$), then
$f$ must be constant.

\section{Metrized graphs versus electrical networks}
\label{Electric_Section}

We now take a moment to give some physical intuition about the Laplacian coming from the
theory of electrical networks. (For a more detailed account of the theory of electrical
networks, see \cite{Bo}, \cite{DS}, and \cite{CR}.) For our purposes, a
\textit{(resistive) electrical network} is a physical model of a metrized graph $\Gamma$
obtained by viewing the vertices of $\Gamma$ as
nodes of the network and the segments of $\Gamma$ as wires, each
with a resistance given by its length.

Using an external device (such as a battery), one can force {\em current} to flow through
the network; for simplicity, we consider only the case where a quantity $I>0$ of current
enters the circuit at some point $a$ and exits at some point $b$.  At all other points of
$\Gamma$, we have {\em Kirchhoff's current law}: The total current flowing into any node
equals the current flowing out of any node.  Mathematically, current is a function which
assigns to each oriented segment $e$ of $\Gamma$ a real number $i_e$, the {\em current
flow} across $e$. {\em Kirchhoff's node law} says that it is possible to define an {\em
electric potential function} $\phi(x) \in \aff$ such that for every oriented segment $e$,
$\phi'_e(x) = -i_e$.  (The minus sign is due to the convention that current flows from
high potential to low potential.) In particular, if $p$ is the initial endpoint and $q$
the terminal endpoint of an oriented segment $e$, then {\em Ohm's law} $\phi(p) - \phi(q)
= i_e L_e$ holds. The potential function $\phi(x)$ is only determined up to an additive
constant; one needs to pick a reference voltage at some point of $\Gamma$ in order to
define the potential at other points.

In our language of directional derivatives, if $p$ is a point of $\Gamma$ and $\vv \in
\vect(p)$ is any direction at $p$, then the current flowing away from $p$ in the
direction $\vv$ is $-D_{\vv} \phi(p)$. Mathematically, Kirchhoff's current law states
that for $p \not\in \{ a,b \}$, we have $-\sigma_p(\phi) = -\sum D_{\vv}\phi(p) = 0$. We
have $-\sigma_a(\phi) > 0$, which says that $a$ is a \textit{current
  source}, and $-\sigma_b(\phi) < 0$, which says that $b$ is a \textit{current
  sink}.  The current entering the network at $a$ is
$-\sigma_a(\phi)$; the current exiting the network at $b$ is $\sigma_b(\phi)$; and we
have $-\sigma_a(\phi) = \sigma_b(\phi) = I$.

Taken together, Kirchhoff's node and potential laws say that given $I>0$, there is a
function $\phi \in \aff$ such that $\Delta \phi = I \cdot \delta_a - I \cdot \delta_b$.
(This will be proved mathematically as a consequence of Corollary~\ref{JFunctionCor} in
\S\ref{jFunctionSection}.)  Note that $\phi$ is determined up to an additive constant by
Theorem~\ref{Laplacian Fact 3}.  Note also that we initially required $-\sigma_a(\phi) =
\sigma_b(\phi) = I$ ({\em conservation of current}), which is demanded mathematically by
Corollary~\ref{Laplacian Fact 2}.

In accordance with physical intuition, the Maximum Principle
(Theorem~\ref{Laplacian Fact 1}) implies that the electric potential
in the network is highest at $a$ (where current enters) and lowest at
$b$ (where it exits).  By convention, one often sets the potential at
$b$ to be zero, in which case we say that the node $b$ is {\em
grounded}.


\section{The Laplacian on a weighted graph} \label{Laplacian Matrix}

In this section, we explain some connections between the classical Laplacian matrix on a
weighted graph and the Laplacian on a metrized graph.

Suppose $G$ is a weighted graph with vertex set $V(G)=\{v_i\}$, edge
set $E(G)=\{e_k\}$, and weights $\{w_{e_k}\}$. If the edge $e_k$ has
endpoints $v_i$ and $v_j$, then we will use the notation
$w_{ij}=w_{e_k}=w_{ji}$ to show the dependence of the weights on the
vertices. For convenience, we set $w_{ij}=0$ if $v_i$ and $v_j$ are not
connected by an edge. In particular $w_{ii}=0$ for all $i$.

\begin{define}
The \textit{Laplacian matrix} associated to a weighted graph $G$ is the $n
\times n$ matrix $Q$ with entries
\[
Q_{ij} = \begin{cases} \sum_k w_{ik}, & \textrm{if $i = j$}
\\ -w_{ij}, & \textrm{if $i \not= j$.} \end{cases} \]
\end{define}
We should note that in the literature, our $Q$ is often called the \textit{combinatorial
Laplacian} or Kirchhoff matrix (see e.g., \cite{Bo}).

The Laplacian matrix encodes interesting information about the graph $G$ (see e.g.,
\cite{Mo}, \cite[\S13]{GR}). For example, zero appears as an eigenvalue of $Q$ with
multiplicity equal to the number of connected components of $G$ (so exactly once in our
case).  Kirchhoff's famous Matrix-Tree Theorem (see \cite[Corollary 13, Chapter II]{Bo})
equates the weighted number of spanning trees of the graph with the absolute value of the
determinant of the matrix obtained by deleting any row and column from $Q$.

Returning to metrized graphs, we've already noted in Exercise~\ref{AffineDetermined}
that a function $f \in \aff$ is completely determined
by its values on the finite set $V_f(\Gamma)$.
Thus, a piecewise affine function on $\Gamma$ yields a function on the
vertices of a certain model for $\Gamma$, and conversely, given a model $G$ and a
function on $V(G)$, we can linearly interpolate to obtain a piecewise
affine function on $\Gamma$.  Our two notions of Laplacian honor this
correspondence:

\begin{thm}\label{affine equivalence}
Suppose $\Gamma$ is a metrized graph, $f \in \aff$, and $G$ is a model of
$\Gamma$ with vertex set $V_f(\Gamma)=\{v_1, \ldots, v_n\}$. Let $\ff$ be the
$n \times 1$ vector with $\ff_i=f(v_i)$. Then
\[
\Delta f = \sum_{i} \left[Q\ff \right]_i \delta_{v_i}.
\]
\end{thm}

\begin{proof}
We already know that $\Delta f$ is discrete if $f$ is piecewise affine. So
it suffices to show that $\left[Q\ff \right]_i= -\sigma_{v_i}(f)$ for any vertex
$v_i$. To that end, we parametrize each segment $e$ having $v_i$ as an
endpoint so that $s_e(0)=v_i$. As $f$ is piecewise affine, the directional
derivatives of $f$ at $v_i$ are given by $f'_e(0) = [f_e(L_e)-f_e(0)]/L_e$.
Recall that the weight of an edge is the reciprocal of its length. We
conclude that
\begin{align*}
\sigma_{v_i}(f) &= \sum_{\substack{\textrm{segments } e \\
\textrm{adjacent to } v_i}} \frac{f_e(L_e) -
f_e(0)}{L_e} =\sum_j w_{ij}\left\{ f(v_j) -f(v_i)\right\} \\
&=-\left\{\left(\sum_k
w_{ik}\right)f(v_i) - \sum_{j \not= i} w_{ij}f(v_j)\right\}  = -\left[Q\ff\right]_i.
\end{align*}
\end{proof}


As a bonus, we now deduce a few useful facts about the Laplacian matrix:
\begin{cor} \label{Bonus Facts}
If $G$ is a weighted graph with $n\times n$ Laplacian matrix $Q$, then
\begin{itemize}
    \item[(i)] The kernel of $Q$ is $1$-dimensional with basis $[1, \ldots, 1]^t$.
    \item[(ii)] If $x \in \RR^n$ is a vector, then $\sum_i \left[Qx\right]_i = 0$.
\end{itemize}
\end{cor}

\begin{proof}
Identify $\RR^n$ with the $n$-dimensional vector space spanned by the
vertices of $G$. A vector $x \in \RR^n$ can be interpreted as a function on
the vertices of $G$, and this function can be linearly interpolated to yield a
piecewise affine function $f$ on the associated metrized graph $\Gamma$. If
$Qx=0$, then Theorem~\ref{affine equivalence} implies that $\Delta f=0$. The Maximum
Principle shows $f$ must be constant, so $x=[c, \ldots, c]^t$ for
some real number $c$. This proves (i). For (ii), use
Corollary~\ref{Laplacian Fact 2} and Theorem~\ref{affine equivalence} to get
\[
\sum_i \left[Qx\right]_i = \int_{\Gamma} \Delta f = 0.
\]
\end{proof}

Now we know the relationship between the Laplacian operator acting on $\aff$
and the Laplacian matrix. In fact, one can prove
that the Laplacian of a piecewise smooth function $f$ is a limit of
Laplacians of piecewise affine approximations of $f$.
To state the result, we introduce the following notation.  If $f
\in \smooth$ and $G_N$ is a model of $\Gamma$ whose vertices contain
$V_f(\Gamma)$, define $f_N$ to be the unique piecewise affine function
with $f_N(p)=f(p)$ for each vertex $p$ of $G_N$ (restrict $f$ to the
vertices of $G_N$ and linearly interpolate).

\begin{thm} \label{Laplacian Convergence}
Suppose $f \in \smooth$. There exists a sequence of models $\{G_N\}$ for
$\Gamma$ such that for all continuous functions $g$ on $\Gamma$, we have
\[
\int_{\Gamma} g \ \Delta f_N \longrightarrow \int_{\Gamma} g \Delta f \qquad \textrm{as $N \to \infty$.}
\]
\end{thm}
That is, the sequence of measures $\{\Delta f_N\}$ converges weakly to $\Delta f$ on
$\Gamma$. By Theorem~\ref{affine equivalence} the discrete measures $\Delta f_N$ can be
computed using the Laplacian matrix.

Theorem~\ref{Laplacian Convergence} is not hard to prove, but we will
not give the proof here.  (A complete proof can be found in \cite{Fa}.)
We mention the theorem in order to display the very close
connection between the Laplacian matrix on a weighted graph and the Laplacian operator on
a metrized graph.


\section{The $j$-function}
\label{jFunctionSection}

In this section, we introduce a three-variable function $j_z(x,y)$ on the metrized graph
$\Gamma$ which allows us, in a sense to be made precise, to invert the Laplacian
operator.\footnote{In comparison with the Riemannian manifold setting, integrating
$j_z(x,y)$ will yield the associated {\em Green's function} for the metrized graph
Laplacian.}
Let $\Meas_0(\Gamma)$ denote the space of measures
of total mass zero on $\Gamma$.
We know from Corollary~\ref{Laplacian Fact 2}
that if $f \in \smooth$ then $\Delta f \in \Meas_0(\Gamma)$.
The following result is a partial converse to this fact.

\begin{thm} \label{DiscreteInversionTheorem}
Let $\nu = \sum c_i \delta_{p_i} \in \emph{Meas}_0(\Gamma)$ be a discrete measure.  Then
there exists a piecewise affine function $f$ on $\Gamma$ such that $\Delta f = \nu$.
\end{thm}

\begin{proof}
Let $S = \{ p_1,\ldots,p_k \}$, and fix a model $G$ for $\Gamma$ with vertex set $V(G)$
containing $S$.  Let $n = \# V(G)$, and let $W$ be the $n$-dimensional real vector space
spanned by the vertices of $G$, which we identify with $\RR^n$. If $Q$ is the Laplacian
matrix associated to $G$, then we know $\ker(Q)$ is $1$-dimensional by
Corollary~\ref{Bonus Facts}(i). The rank-nullity theorem implies that $\im(Q)$ is
$(n-1)$-dimensional.

By Theorem~\ref{affine equivalence}, solving $\Delta f = \nu$ is equivalent
to finding a vector $x \in W$ with $Qx = [c_1, \ldots, c_n]^t$. Let $W_0$ be
the $(n-1)$-dimensional subspace of $W$ consisting of vectors
$[a_1,\ldots,a_n]^t$ such that $\sum a_i = 0$. Corollary~\ref{Bonus Facts}(ii)
shows $\im(Q)$ is contained in $W_0$. As these two spaces have the same
dimension, they must be equal. The condition $\nu \in \Meas_0(\Gamma)$ says
$\sum c_i = 0$, so $[c_1, \ldots, c_n]^t$ lies in the image of $Q$.
\end{proof}

We now single out a special case of this result which is of particular interest. In what
follows, we write $\Delta_x$ instead of $\Delta$ if we wish to emphasize that we are
taking the Laplacian with respect to the variable $x$.

\begin{cor}
\label{JFunctionCor}
For fixed $y,z \in \Gamma$, there exists a unique piecewise affine function $j(x) =
j_z(x,y)$ satisfying
\[
\Delta_x j_z(x,y) = \delta_y(x) - \delta_z(x), \qquad j_z(z,y) = 0.
\]
\end{cor}

\begin{proof}
The existence of $j(x)$ follows from Theorem~\ref{DiscreteInversionTheorem}, and
uniqueness follows from Theorem~\ref{Laplacian Fact 3}.
\end{proof}

We now justify our assertion that the $j$-function allows us to ``invert the Laplacian''
on the space $\Meas_0(\Gamma)$. Recall from Theorem~\ref{DiscreteInversionTheorem} that
given a discrete measure $\nu \in \Meas_0(\Gamma)$, there exists a function $f \in \aff$
(unique up to an additive constant) that satisfies the differential equation $\Delta f =
\nu$. The next result shows that we can explicitly describe such a function $f$ using the
$j$-function:

\begin{thm}
\label{DiscreteInversionTheorem2} Let $\nu = \sum c_i \delta_{p_i} \in
\emph{Meas}_0(\Gamma)$ be a discrete measure.  Then for any fixed $z
\in \Gamma$,
the function
\[
f(x) = \int_\Gamma j_z(x,y) d\nu(y) = \sum_i c_i j_z(x,p_i)
\]
is piecewise affine and satisfies the equation $\Delta f = \nu$.
\end{thm}

\begin{proof}
The condition $\nu \in \Meas_0(\Gamma)$ means that $\sum c_i = 0$.
Therefore
\[
\Delta f= \sum_i c_i \left(
\delta_{p_i} - \delta_z \right) = \nu.
\]
Since the $j$-function is piecewise affine, $f$ is as well.
\end{proof}

We mention (see \cite{BR} for a proof) that Theorem~\ref{DiscreteInversionTheorem2}
admits the following generalization to arbitrary (not necessarily discrete) measures $\nu
\in \Meas_0(\Gamma)$: For fixed $z \in \Gamma$, the function $f(x) = \int_\Gamma j_z(x,y)
d\nu(y)$ is in $\smooth$ and satisfies the equation $\Delta f = \nu$. In particular, if
$\nu$ is a measure on $\Gamma$, then we can solve the differential equation $\Delta f =
\nu$ {\em if and only if} $\nu \in \Meas_0(\Gamma)$.

The function $j_z(x,y)$ has an interpretation in terms of electrical
networks.  Recalling our description of the electrical network
associated to a metrized graph given in \S\ref{Electric_Section}, the
function $j_z(x,y)$ is the electric potential at $x$ if one unit of
current enters the network at $y$ and exits at $z$, and the node $z$
is grounded.  So one could build a real-life model of the metrized
graph $\Gamma$ with wires, hook up a battery, and empirically
determine the values of the $j$-function!

\begin{exercise}
\label{PositivityExercise}
Physical intuition suggests that the $j$-function should be
nonnegative; prove more precisely that
\[
0 \leq j_z(x,y) \leq j_z(y,y)
\]
for all $x,y,z \in \Gamma$.
[{\bf Hint:} For fixed $y$ and $z$, apply the Maximum Principle to $j_z(x,y)$ and its
negative.]
\end{exercise}

The three-variable function $j_z(x,y)$ satisfies a magical four-term
identity, which will be used in various guises throughout this
section and the next.  The proof of this identity
is an excellent illustration of the theory developed in \S\ref{Continuous Laplacian}.

\begin{thm}[Magical Identity]
For all $x,y,z,w \in \Gamma$, we have the identity
\begin{equation*}
\label{MagicalIdentity}
j_z(x,y) -j_z(w,y)= j_w(y,x) - j_w(z,x).
\end{equation*}
\end{thm}

\begin{proof}
Fix $x,y,z,w \in \Gamma$. On one hand, we have
\[
\int_{\Gamma} j_z(u,y) \Delta_u \left( j_w(u,x) \right) = \int_{\Gamma} j_z(u,y)
\left\{\delta_x(u)-\delta_w(u)\right\} = j_z(x,y) - j_z(w,y).
\]
By Theorem~\ref{Laplacian Fact 1}, this is equal to
\[
\int_{\Gamma} j_w(u,x) \Delta_u \left( j_z(u,y)\right)  = \int_{\Gamma}
 j_w(u,x) \left\{\delta_y(u)-\delta_z(u)\right\} = j_w(y,x)-j_w(z,x).
 \]
\end{proof}

The Magical Identity allows us to prove two useful symmetries for the
$j$-function.

\begin{cor}
\label{JFunctionSymmetries} For $x,y,z \in \Gamma$, the $j$-function satisfies
\begin{itemize}
\item[(i)] $j_z(x,y) = j_z(y,x)$
\item[(ii)] $j_z(x,x) = j_x(z,z)$
\end{itemize}
\end{cor}


\begin{proof}
For (i), if we set $w=z$ in the Magical Identity, we obtain
\[
j_z(x,y) -j_z(z,y) = j_z(y,x) - j_z(z,x).
\]
Since $j_z(z,x) =j_z(z,y)= 0$, the result follows.

For (ii), substitute $x=z, y=w$ into the Magical Identity to get
\[
j_z(z,w) -j_z(w,w)= j_w(w,z) - j_w(z,z).
\]
Since $j_z(z,w) = j_w(w,z) =  0$, the result follows by swapping $w$ for $x$.
\end{proof}


In passing, we mention that $j_z(x,y)$ has a very strong continuity property: it is
jointly continuous in $x,y,$ and $z$. That is, the value of the $j$-function varies
continuously if we make small variations to $x,y,$ and $z$ simultaneously. Our electrical
network interpretation makes this statement quite plausible:
the value on our
voltmeter should vary continuously when we move the battery terminals
and the point at which we're reading the voltage.
A mathematical proof is outlined in the next exercise (see \cite{CR}
for a different approach).


\begin{exercise}
\label{JointContinuityExercise}
\begin{itemize}  \item[] \hspace{1in}
\item[(a)] Let $I,I'$ be closed intervals in $\RR$. Suppose $f : I
\times I' \to \RR$ has the property that $f(x,y)$ is affine in $x$ and $y$ separately.
Then $f(x,y) = c_1 + c_2 x + c_3 y + c_4 xy$ for some $c_1,\ldots,c_4 \in \RR$.
\item[(b)] Use (a) to show that for fixed $z \in \Gamma$, $j_z(x,y)$
  is jointly continuous as a function of $x$ and $y$.
\item[(c)] Use Theorem~\ref{Laplacian Fact 3} to prove the five-term identity
\[
j_z(x,y) = j_w(x,y) -j_w(x,z) - j_w(z,y) + j_w(z,z).
\]
\item[(d)] Deduce from (b) and (c) that $j_z(x,y)$ is jointly
  continuous in $x,y$, and $z$.
\end{itemize}
\end{exercise}


We now define another useful function motivated by the theory of
electrical networks:
\begin{define}
The \textit{effective resistance} between two points $x,y$ of a metrized graph is given
by
\[
r(x,y)=j_y(x,x)=j_x(y,y).
\]
\end{define}

The fact that $j_y(x,x)=j_x(y,y)$ is just a restatement of the second
symmetry of the $j$-function in Corollary~\ref{JFunctionSymmetries}. In terms of
electrical networks, the effective resistance between two nodes $x$ and
$y$ is the absolute value of the potential difference between $x$ and $y$ when a unit
current enters the network at $x$ and exits at $y$.

We now introduce some useful techniques for calculating the $j$-function and
the effective resistance function. Rules (ii) and (iii) in Theorem~\ref{CircuitTheoryProp} below are
essentially the familiar series and parallel transforms from circuit
theory. The proofs of Theorems \ref{CircuitTheoryProp} and
\ref{ExplicitResistanceTheorem} below are adapted from \cite{Zhang}.

A \textit{subgraph} of the metrized graph $\Gamma$ is a subspace of $\Gamma$ which is a metrized
graph in its own right. In the statement of Proposition~\ref{CircuitTheoryProp},
$\Gamma_1$ and $\Gamma_2$ will always denote subgraphs of $\Gamma$.  We let $j_z(x,y)$
(resp. $j_{z,1}(x,y), j_{z,2}(x,y)$) denote the $j$-function on $\Gamma$ (resp. on
$\Gamma_1,\Gamma_2$), and similarly we let $r(x,y)$ (resp. $r_1(x,y), r_2(x,y)$) denote
the effective resistance function on $\Gamma$ (resp. on $\Gamma_1,\Gamma_2$).

\begin{figure}[!ht]
  \begin{picture}(200,60)(0,0)
    \put(-50,10){\scalebox{0.6}{\includegraphics{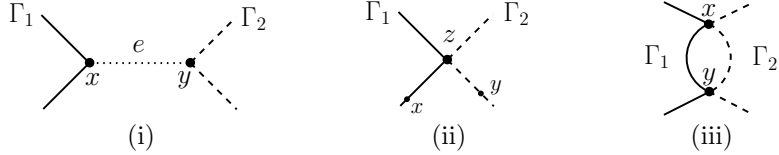}}}

    \put(-5,0){\textrm{(i)}}
    \put(110,0){\textrm{(ii)}}
    \put(208,0){\textrm{(iii)}}

  \end{picture}

  \caption{
These three figures illustrate the three parts of Proposition~\ref{CircuitTheoryProp}.
      The solid lines (resp. dashed lines) indicate the segments of the diagram belonging to $\Gamma_1$
      (resp. to $\Gamma_2$).
}
  \label{Circuit Laws}
\end{figure}

\begin{thm}
\label{CircuitTheoryProp}
Let $\Gamma$ be a metrized graph, and let $\Gamma_1$ and $\Gamma_2$ be
subgraphs.
\begin{itemize}
\item[(i)] Suppose $e$ is a segment in $\Gamma$ of length $L$ with
  endpoints $x,y$, and assume that $\Gamma = \Gamma_1 \cup \Gamma_2
  \cup e$ with $\Gamma_1 \cap e = \{ x \}$, $\Gamma_2 \cap e = \{ y
  \}$, and $\Gamma_1 \cap \Gamma_2 = \emptyset$.
(Compare Figure~\ref{Circuit Laws}(i).)
Then $r(x,y) = L$.
\item[(ii)] Suppose $\Gamma = \Gamma_1 \cup \Gamma_2$ with $\Gamma_1
  \cap \Gamma_2 = \{ z \}$.
(Compare Figure~\ref{Circuit Laws}(ii).)
Then for all $x \in \Gamma_1$ and $y \in
  \Gamma_2$, we have $r(x,y) = r_1(x,z) + r_2(z,y)$.
\item[(iii)] Suppose $\Gamma = \Gamma_1 \cup \Gamma_2$ with $\Gamma_1
  \cap \Gamma_2 = \{ x,y \}$.
(Compare Figure~\ref{Circuit Laws}(iii).)
Then
\[
\frac{1}{r(x,y)} = \frac{1}{r_1(x,y)} + \frac{1}{r_2(x,y)}.
\]
\end{itemize}
\end{thm}

\begin{proof}

For (i), we pick a parametrization $s_e:[0,L] \to e$ such that $s_e(0)=x$ and $s_e(L)=y$. Let
$t:e \to [0,L]$ be the inverse of $s_e$.  We claim that
\[
j_x(z,y) = \begin{cases}
    0, & \textrm{if }z \in \Gamma_1 \\
    t(z), & \textrm{if }z \in e \\
    L, & \textrm{if }z \in \Gamma_2.
\end{cases}
\]

Indeed, it is easily verified that the Laplacian of the right-hand side with respect to
$z$ is $\delta_y - \delta_x$, and that the two sides agree when $z = x$.  The claim
therefore follows from Theorem~\ref{Laplacian Fact 3}, and the desired result follows by
setting $z = y$.

For (ii), we claim that
\[
j_x(w,y) = \begin{cases}
    j_{x,1}(w,z), & \textrm{if }w \in \Gamma_1 \\
    r_1(x,z) + j_{z,2}(w,y), & \textrm{if }w \in \Gamma_2. \end{cases}
\]
The point is that the right-hand side is continuous at $w = z$, has Laplacian
 equal to $(\delta_z - \delta_x) + (\delta_y - \delta_z) = \delta_y -
 \delta_x$, and is zero when $w = x$.
The result then follows by setting $w = y$.

We leave the proof of (iii) as an exercise for the reader.

\end{proof}

\begin{exercise}
Verify part (iii) of Theorem~\ref{CircuitTheoryProp} by first showing that
\begin{large}
\[
j_x(z,y) = \begin{cases}
    \frac{r_2(x,y)}{r_1(x,y) + r_2(x,y)} j_{x,1}(z,y), &\textrm{if } z \in \Gamma_1 \\
    \hspace{.5in} & \hspace{.5in} \\
    \frac{r_1(x,y)}{r_1(x,y) + r_2(x,y)} j_{x,2}(z,y), &\textrm{if } z \in \Gamma_2.
\end{cases}
\]
    \end{large}
\end{exercise}

\begin{exercise}
Show that the function $r(x,y)$ is jointly continuous in $x$ and $y$,
and that for fixed $y \in \Gamma$, $r(x,y)$ is continuous and piecewise
quadratic in $x$.
[{\bf Hint:} Use Exercise~\ref{JointContinuityExercise}.]
\end{exercise}

Using Theorem~\ref{CircuitTheoryProp}, we can derive an explicit description of the
function $r(x,y)$ when $x$ varies along a single segment of $\Gamma$ having $y$ as an
endpoint. To state the result, we define a quantity $R_e$ associated to a segment $e$ of
$\Gamma$ as follows. Let $\inter{e}$ denote the interior of the segment $e$, and let
$\Gamma_e$ be the complement of $\inter{e}$ in $\Gamma$. If $\Gamma_e$ is connected, then
$\Gamma_e$ is a subgraph of $\Gamma$, and we define $R_e$ to be the effective resistance
$r(y,z)$ between the endpoints $y$ and $z$ of $e$ computed on $\Gamma_e$. If $\Gamma_e$
is not connected (i.e., if $e$ is not part of a cycle), we define $R_e$ to be $\infty$.
Loosely speaking, $R_e$ is the effective resistance between the endpoints of $e$ in the
subgraph obtained by deleting $e$.

The next result is motivated by the following intuition: to calculate
$r(x,y)$ on the segment $e$, we can think of $e$ and its complement $\Gamma_e$
as being connected in parallel, and $x$ splits $e$ into two segments
connected in series. We can then use the parallel and series transforms to
calculate $r(x,y)$.

\begin{thm}
\label{ExplicitResistanceTheorem} Let $e$ be a closed segment of $\Gamma$ of length
$L_e$, let $y,z$ be the endpoints of $e$, and parametrize $e$ by $s_e:[0,L_e]\to e$ with
$s_e(0)=y$ and $s_e(L_e)=z$. Suppose $t:e \to [0,L_e]$ is the inverse of $s_e$. Then for $x \in
e$, we have
\[
r(x,y) = t(x) - \frac{1}{L_e + R_e} t(x)^2,
\]
where $\frac{1}{L_e + R_e} = 0$ if $R_e = \infty$.
\end{thm}

\begin{proof}
If $\Gamma_e$ is not connected, then $\Gamma = \Gamma_1 \cup \Gamma_2 \cup e$, with
$\Gamma_1 \cap e = \{ y \}$, $\Gamma_2 \cap e = \{ z \}$, and $\Gamma_1 \cap \Gamma_2 =
\emptyset$. Then $r(x,y) = t(x)$ by part (i) of Proposition~\ref{CircuitTheoryProp}.

Now suppose that $\Gamma_e$ is connected. Then $x$ breaks $e$ into two closed segments
$\Gamma_1 = [t(y),t(x)]$ and $\Gamma_2 = [t(x),t(z)]$, and $\Gamma = \Gamma_1 \cup
\Gamma_2 \cup \Gamma_e$.  Letting $\Gamma_3 = \Gamma_2 \cup \Gamma_e$, we have (with the
obvious notation):
\begin{equation*}
\begin{aligned}
\frac{1}{r(x,y)} &= \frac{1}{r_1(x,y)} + \frac{1}{r_3(x,y)},
\qquad\text{by Prop.~\ref{CircuitTheoryProp}(iii)} \\
&= \frac{1}{r_1(x,y)} + \frac{1}{r_2(x,z) + R_e},
\qquad\text{by Prop.~\ref{CircuitTheoryProp}(ii)} \\
&= \frac{1}{t(x)} + \frac{1}{L_e - t(x) + R_e},
\qquad\text{by Prop.~\ref{CircuitTheoryProp}(i)}. \\
\end{aligned}
\end{equation*}

The desired formula now follows because of the simplification
\[
\left( \frac{1}{t(x)} + \frac{1}{L_e - t(x) + R_e} \right)^{-1} = t(x) - \frac{1}{L_e +
R_e} t(x)^2.
\]
\end{proof}

\begin{exercise}
\label{CurrentExercise}
Let $\Gamma$ be a metrized graph of total length $L$.
Fix $a,b \in \Gamma$, and choose a
vertex set $V(\Gamma)$ containing $a$ and $b$.
Let $e$ be an oriented segment of $\Gamma$
beginning at $p$ and ending at $q$.
For $x \in \Gamma$, let $\phi(x) = j_b(x,a)$, and
define $i_e = \frac{\phi(p) - \phi(q)}{L_e} = -\phi_e' $.
(In terms of electrical networks, $i_e$ is the current flowing across $e$ when
a unit current enters the network at $a$ and exits at $b$.)  Also, define $r(e) = r(p,q)$.
\begin{itemize}
\item[(a)] Show that $r(e) \leq L_e$. [{\bf Hint:} Use Prop.~\ref{CircuitTheoryProp}.]
\item[(b)] Show that $r(x,y)$ is a metric on $\Gamma$.
   [{\bf Hint:} For the triangle inequality, use Exercise~\ref{JointContinuityExercise}.]
\item[(c)] Deduce that $r(a,b)$ is bounded above by the length of any path from $a$ to $b$,
   and conclude that $0 \leq r(x,y) \leq L$ for all $x,y \in \Gamma$.
\end{itemize}
\end{exercise}


\section{The canonical measure and Foster's Theorem} \label{Canonical Measure Section}

Calculating the Laplacian of the effective resistance function $r(x,y)$
for fixed $y$ is not so easy just from the definitions, but our explicit description in
Theorems~\ref{ExplicitResistanceTheorem}
and \ref{IndependenceTheorem} below will allow us to do it in a slick way.
The first half of this section will be devoted to figuring out
$\Delta_x r(x,y)$, and in the second
half we reap the benefits of this calculation by proving some
interesting results from graph theory, including Foster's theorem.
The method presented here is a simplified version of \S2 of \cite{CR}.

\begin{example}
If $\Gamma = [0,1]$, then Theorem~\ref{CircuitTheoryProp}(i) shows
that $r(x,y) = |x-y|$, and a
simple calculation shows that $\Delta_x r(x,y) = \delta_0(x) + \delta_1(x) -
2\delta_y(x)$. Interestingly, we see that $\Delta_x r(x,y) + 2\delta_y(x)$ is independent
of $y$. This simple example actually illustrates a general phenomenon.
\end{example}

\begin{thm}
\label{IndependenceTheorem}
For any metrized graph, $\Delta_x r(x,y) + 2\delta_y(x)$ is a measure which
is independent of $y$.
\end{thm}

\begin{proof}
Let $z,w \in \Gamma$ be arbitrary. Set $x=y$ in the Magical Identity of
\S\ref{jFunctionSection} to get
\[
j_z(y,y)-j_z(w,y)=j_w(y,y)-j_w(z,y).
\]
Applying Corollary~\ref{JFunctionSymmetries}, we obtain
\[
r(y,z)-j_z(y,w)=r(y,w)-j_w(y,z).
\]
Taking the Laplacian of both sides with respect to $y$ and recalling that $\Delta_y
j_z(y,x) = \delta_x - \delta_z$, we get
\[
\Delta_y r(y,z) - \delta_w+\delta_z = \Delta_y r(y,w) - \delta_z +\delta_w.
\]
Rearranging, we see that $\Delta_y r(y,z) + 2\delta_z = \Delta_y r(y,w)
+ 2 \delta_w$.
As $w$ and $z$ were arbitrary, the result follows.
\end{proof}

\begin{define} \label{canonicalmeasure}
The \textit{canonical measure} on a metrized graph $\Gamma$
is given by
\[
\mucan = \frac{1}{2}\Delta_x r(x,y) + \delta_y (x),
\]
where $y \in \Gamma$ is arbitrary.
Theorem~\ref{IndependenceTheorem} shows that $\mucan$ is independent of the choice of $y$.
\end{define}

Recall from Corollary~\ref{Laplacian Fact 2} that $\Delta_x r(x,y)$ is a measure of total
mass zero, so we see from Definition~\ref{canonicalmeasure}
that $\mucan$ has total mass~$1$.

\begin{exercise}
\label{TauExercise}
Show that the quantity
\[
\tau(\Gamma) = \frac{1}{2} \int_\Gamma r(x,y) d\mucan(x)
\]
is independent of the choice of $y \in \Gamma$.
[{\bf Hint:} Use Theorem~\ref{Laplacian Fact 1}.]
\end{exercise}

We now give an explicit description of the measure $\mucan$.

\begin{thm}
\label{CanonicalMeasureTheorem}
Let $n_p$ denote the valence of a vertex $p \in V(\Gamma)$.  Then
\[
\mucan = \sum_{\textrm{\emph{vertices} } p} \left( 1-\frac{1}{2}n_p \right) \delta_p +
\sum_{\textrm{\emph{segments} } e} \frac{1}{R_e + L_e} dx\vert_{e}.
\]
\end{thm}

\begin{proof}We compute the discrete and continuous parts of $\mucan$ separately.

\textit{Continuous part:}
Let $e$ be an oriented segment of $\Gamma$ which begins at $y$ and
ends at $z$.
If $x$ lies on $e$, then we're in the situation of
Theorem~\ref{ExplicitResistanceTheorem}, and we calculate that
\begin{equation*} \tag{$\dagger$} \label{eq:explicitresistance}
\Delta_x \left\{r(x,y)\vert_e \right\} = \frac{2}{L_e+R_e} dx\vert_e + \delta_z -
\delta_y.
\end{equation*}
Since $\mucan = \frac{1}{2} \Delta_x r(x,y) + \delta_y$ is independent of our choice of
$y$, \eqref{eq:explicitresistance} shows that the continuous part of $\mucan$ along $e$
must be $\frac{1}{L_e+R_e} dx\vert_e$.

\textit{Discrete part:}
If $y$ is an endpoint of a segment $e$, then $r(x,y)$ is quadratic
along the interior of $e$ by Theorem~\ref{ExplicitResistanceTheorem}.
It follows that the discrete part of $\mucan$ is supported on
$V(\Gamma)$.  Let $p \in V(\Gamma)$ be a vertex.
Using (\ref{AlternateDefinition}) in \S\ref{Continuous Laplacian},
we calculate from Equation \eqref{eq:explicitresistance}
that $\frac{1}{2} \Delta_x r(x,p)$ contributes $-\frac{1}{2} \delta_p$ to the discrete part of
$\mucan$ at $p$ for each segment $e$ beginning at $p$. Recalling that $\mucan=\frac{1}{2}
\Delta_x r(x,p) + \delta_p$, the coefficient of $\delta_p$ in $\mucan$
must therefore be $1 - \frac{1}{2} n_p.$

\end{proof}

\begin{example}
If $\Gamma$ is a circle of length 1, then every vertex has
valence~$2$, so $\mucan$ has no discrete part. For the continuous part,
divide the circle into three segments $e_1, e_2, e_3$, each of length $1/3$.
Then we get \[ \mucan = dx\mid_{e_1} + dx\mid_{e_2} + dx \mid_{e_3} = dx. \]
\end{example}

\begin{example}
Let $\Gamma$ be the star of Figure~\ref{Flux Laplacian Example}. Then
$\mucan$ has no continuous part because $R_e$ is infinite for all edges.
Therefore \[ \mucan =
\frac{1}{2}\delta_{P}-\frac{1}{2}\delta_{Q}+\frac{1}{2}\delta_{R}+\frac{1}{2}\delta_{S}.
\]
\end{example}

Theorem~\ref{CanonicalMeasureTheorem} has some interesting consequences for weighted
graphs.  For example, we have the following result from \cite{CR}:

\begin{cor}
\label{CanonicalMeasureCor} Let $G$ be a weighted graph with vertex set $V(G)$ and edge
set $E(G)$. Then
\[
\sum_{\textrm{\emph{edges} } e} \frac{L_e}{R_e + L_e} = 1+ \#E(G) - \#V(G).
\]
\end{cor}

\begin{proof}
Integrating both sides of the formula in Theorem~\ref{CanonicalMeasureTheorem} over
$\Gamma$, we obtain:
\[
1 = \sum_{\textrm{vertices } p} \left( 1 - \frac{1}{2} n_p \right) + \sum_{\textrm{edges
} e} \frac{L_e}{R_e + L_e}.
\]
(Here we have summed over edges of $G$ instead of segments of $\Gamma$, but the two sets
are in bijective correspondence.) As each edge in $G$ connects exactly 2 vertices,
we have
\[
\sum_{\textrm{vertices } p} n_p = 2 \left\{\# E(G)\right\}.
\]
Therefore
\[
1 = \# V(G) - \# E(G) + \sum_{\textrm{edges } e} \frac{L_e}{R_e + L_e},
\]
which is equivalent to the desired formula.
\end{proof}

It is a well-known fact from graph theory that $1 + \# E(G) - \# V(G)$
is the number of linearly independent cycles on $G$ (see \cite[Theorem
9, Chapter II]{Bo}). This is a topological invariant which only depends on the associated metrized
graph $\Gamma$.

\begin{cor}[Foster's Theorem]
\label{Foster} For an edge $e$ in a weighted graph $G$, let
$r(e)$ denote the effective resistance $r(x,y)$ between the endpoints
$x$ and $y$ of $e$ on the associated metrized graph $\Gamma$.
Let $w_e = 1/L_e$ be the weight of the edge $e$.
Then
\[
\sum_{\textrm{\emph{edges} }e} w_e \ r(e) =
\sum_{\textrm{\emph{edges} }e} \frac{r(e)}{L_e} = \#V(G) - 1.
\]
\end{cor}

\begin{proof}
If $R_e = \infty$, then $r(e)= L_e$ by Proposition~\ref{CircuitTheoryProp}(i).
Otherwise, by Theorem~\ref{ExplicitResistanceTheorem}, we have
\[
r(e) = L_e-\frac{L_e^2}{L_e+R_e}=\frac{L_e R_e}{L_e+R_e}.
\]
Combining these observations, we see that
\begin{align*}
\sum_{\textrm{edges } e} \frac{r(e)}{L_e} &=
\sum_{\substack{\textrm{edges } e \\ \textrm{with } R_e = \infty }} 1
+ \sum_{\substack{\textrm{edges } e \\ \textrm{with } R_{e} \neq \infty}} \frac{R_e}{L_e+R_e}
\\
&= \#E(G) + \sum_{\substack{\textrm{edges } e \\ \textrm{with } R_e \neq \infty }}
\left\{\frac{R_e}{L_e + R_e}-1\right\}\\
&= \#E(G) - \sum_{\textrm{edges } e} \frac{L_e}{L_e+R_e}.
\end{align*}
The result follows immediately from Corollary~\ref{CanonicalMeasureCor}.
\end{proof}

\begin{example}
If $G$ is a tree, then we have
$r(e) = L_e$ for all $e$ by Proposition~\ref{CircuitTheoryProp}(i),
and $\# E(G) = \# V(G) - 1$.  Therefore
$\sum_e \frac{r(e)}{L_e} = \# E(G) = \# V(G) - 1$ as predicted by
Foster's theorem.

More generally, for arbitrary $G$ it follows from Exercise~\ref{CurrentExercise}(a) that
$0 \leq \frac{r(e)}{L_e} \leq 1$ for each edge $e$,
so that \textit{a priori} we have $\sum_e \frac{r(e)}{L_e} \leq \# E(G)$.
Foster's theorem is equivalent to the assertion that the difference
$\# E(G) - \sum_e \frac{r(e)}{L_e}$ is equal to the number of
independent cycles in $G$.
\end{example}

\begin{example}
Foster's theorem can be a useful tool for calculating effective
resistances, especially in the presence of symmetry.  For example, let $G = K_n$ be the complete graph
on $n\geq 2$ vertices, with all edge weights equal to 1.  By symmetry,
the effective resistance $r(x,y)$ between distinct points $x,y \in V(G)$ is
independent of $x$ and $y$; let $r$ denote the common value.  Foster's theorem gives
\[
\sum_{\textrm{edges }e} w_e \ r(e) = \binom{n}{2} \cdot r = n - 1,
\]
so that $r = 2/n$.
\end{example}

\section{Eigenfunctions of the Laplacian} \label{Fourier Analysis}

Suppose $\Gamma$ is a circle of length 1.  Then for $f \in \smooth$, we have \[\Delta f =
-f''(x) dx + (\textit{discrete measure}).\] A standard computation shows that the nonzero
piecewise smooth functions $\phi$ that satisfy the equation
\[
\Delta \phi = \lambda \phi(x) dx
\]
for some $\lambda \in \RR$ are precisely the constant multiples of $\sin(2\pi nx)$ and
$\cos(2\pi nx)$ for $n \in \ZZ$. These functions will be called the
\textit{eigenfunctions of the Laplacian} on $\Gamma$. The corresponding
\textit{eigenvalues} are $\lambda_n = 4\pi^2 n^2$.

It is convenient to normalize each eigenfunction $\phi$ of the Laplacian so that
$\int_\Gamma \phi(x)^2 \, dx = 1$; i.e., $\phi$ has $L^2$-norm equal to~1. By standard
calculus facts, for $m \not= n$ and $n \not= 0$ we have
\begin{align*}
\int_\Gamma \sin^2(2\pi nx) \, dx &= \int_\Gamma \cos^2(2\pi nx) \, dx
= \frac{1}{2}, \\
\int_\Gamma \sin(2\pi nx)\sin(2\pi mx) \, dx &= \int_\Gamma \cos(2\pi
nx)\cos(2\pi mx) \, dx = 0,\\
\int_\Gamma \sin(2\pi nx)\cos(2\pi mx) \, dx &= 0. \\
\end{align*}
Therefore
\[
\Lambda = \{1\} \cup \{ \sqrt{2}\cos(2\pi nx) \}_{n \geq 1} \cup \{ \sqrt{2}\sin(2\pi nx) \}_{n
\geq 1}
\]
is an $L^2$-\textit{orthonormal} set of eigenfunctions for the Laplacian.

We make the following observations about the set $\Lambda$:
\begin{itemize}
\item Each $L^2$-normalized eigenfunction of $\Delta$ occurs exactly once on this
list.
\item Each nonzero eigenvalue $4\pi^2 n^2$ occurs twice, and
the eigenvalue~$0$ occurs with multiplicity~$1$.
\end{itemize}

A standard result in Fourier analysis is the following:

\begin{thm}
\label{StandardFourierTheorem} Let $\Gamma$ be a circle of length 1.
Then any $f \in \smooth$ can be expanded as a uniformly convergent series
\[
f(x) = a_0 + \sum_{n\geq 1} a_n \sqrt{2}\cos(2\pi nx) + \sum_{n\geq 1} b_n
\sqrt{2}\sin(2\pi nx),
\]
where the \textit{Fourier coefficients} $a_n, b_n$ are determined by
\[
\begin{cases}
a_0 = \int_\Gamma f(x) dx & \\
a_n = \int_\Gamma f(x) \sqrt{2}\cos(2\pi n x) dx & n \geq 1\\
b_n = \int_\Gamma f(x) \sqrt{2}\sin(2 \pi n x) dx & n \geq 1.\\
\end{cases}
\]

\end{thm}

Viewed in this light, Fourier analysis on the circle is the theory of eigenfunctions of
the Laplacian on the underlying metrized graph.
\footnote{Of course, there are many variants and generalizations of
Theorem~\ref{StandardFourierTheorem}, and much sophisticated mathematics has been
developed to address what happens if $f$ satisfies hypotheses weaker than piecewise
smoothness (for example, if $f$ is merely continuous).  But this article is not the place to
discuss such matters!}

A nice fact is that one can generalize Fourier analysis to an arbitrary metrized graph.
One way to do this is as follows. Fix a measure $\mu$ on $\Gamma$; for simplicity, we
will assume that $\mu$ has total mass 1. Let
\[
\smoothmu = \{ f \in \smooth \; : \; \int_\Gamma f(x) d\mu(x) = 0 \}.
\]

We now make the following somewhat non-intuitive definition:

\begin{define} \label{Eigenfunction} A nonzero function $\phi \in \smooth$ is an
\textit{eigenfunction of the Laplacian with respect to $\mu$} if $\phi \in \smoothmu$
  and satisfies the equation
\[
\Delta \phi = \lambda \phi(x) dx - C \mu
\]
for some $\lambda, C \in \RR$.
\end{define}

Note that the value of the constant $C$ is completely determined by $\lambda$ and $\phi$
in the above equation.  Indeed, integrating both sides and recalling that $\Delta\phi$
has total mass zero and $\mu$ has total mass one shows that
\[C= \lambda \int_{\Gamma} \phi(x) dx.\]

A sequence $\{ \phi_n \}$ of distinct eigenfunctions is \textit{orthonormal} if
$\int_{\Gamma} \phi_i(x) \phi_j(x)dx = \delta_{ij}$, where $\delta_{ij}=1$ if $i=j$ and
$0$ otherwise. For our purposes, the sequence is called \textit{complete} if every
eigenfunction of the Laplacian with respect to $\mu$ is a scalar multiple of some
$\phi_n$. \footnote{This is a non-standard definition of complete. In \cite{BR} it is
proved that $L^2(\Gamma)$ admits a complete orthonormal basis of eigenfunctions of the
Laplacian in the standard sense of complete; i.e., the only $L^2$-function orthogonal to
all of the eigenfunctions is the zero function.}

It may not be clear \textit{a priori} what role the $C \mu$ term is playing in the
definition of an eigenfunction. However, this definition turns out to be quite flexible and
useful, as illustrated by the following result (see \cite{BR} for a proof):

\begin{thm}
\label{GeneralizedFourierTheorem} Suppose $\Gamma$ is a metrized graph. Let
$\mu$ be a measure of total mass~1 on $\Gamma$, and consider a
complete orthonormal sequence $\{ \phi_n \}_{n \geq 1}$ of eigenfunctions of
the Laplacian with respect to $\mu$.
The corresponding eigenvalues $\{\lambda_n\}$ are all positive and each occurs with finite multiplicity.
Furthermore, every $f \in \smoothmu$ can be expanded as a uniformly convergent
series
\[
f(x) = \sum_{n\geq 1} a_n \phi_n(x),
\]
where the generalized Fourier coefficients~$a_n$ are determined by the formula
\[
a_n = \int_\Gamma f(x) \phi_n(x) dx.
\]
\end{thm}

If we add in the constant function~$\id$, then it follows from
Theorem~\ref{GeneralizedFourierTheorem} that every $f \in \smooth$ can be uniquely
expressed as
\[
f(x) = a_0 + \sum_{n\geq 1} a_n \phi_n(x),
\]
where $a_0 = \int_\Gamma f(x) d\mu(x)$ and the $a_n$'s are as before.

Though the main interest of this result is the fact that it applies to arbitrary metrized
graphs, we illustrate what's happening in Theorem~\ref{GeneralizedFourierTheorem} by
considering the special case where $\Gamma = [0,1]$ and $\mu = \delta_0$ is a point mass
at~$0$.

What are the eigenfunctions of the Laplacian in this case? By
Definition~\ref{Eigenfunction} and the fact that an eigenfunction is required to be in
$\smoothmu$, we demand that
\[
-\phi''(x) dx + \phi'(1)\delta_1 - \phi'(0) \delta_0 = \lambda \phi(x) dx - C \delta_0,
\quad \phi(0) = 0,
\]
for some $\lambda, C \in \RR$.

Thus $\phi''(x) = -\lambda \phi(x)$, $\phi'(1) = 0$, and $\phi(0) = 0$. A computation now
shows that $\phi(x)$ must be a constant multiple of $\sin(\pi nx/2)$ for some odd
positive integer~$n$. It is then easy to verify that
\[
\Lambda = \left\{ \sqrt{2} \sin\left(\frac{\pi nx}{2}\right) \right\}_{n\geq 1 \textrm{ odd}}
\]
forms a complete orthonormal set of eigenfunctions for the Laplacian with
respect to $\delta_0$.  The corresponding eigenvalues are $\pi^2 n^2/4$, each of which occurs
with multiplicity one.

The next result follows immediately from Theorem~\ref{GeneralizedFourierTheorem}, but in
order to show the connection to classical Fourier analysis, we will deduce it directly
from Theorem~\ref{StandardFourierTheorem}.

\begin{thm}
\label{SpecificFourierTheorem} Every $f \in {\emph S}([0,1])$ can be written as a
uniformly convergent generalized Fourier series of the form
\[
f(x) = f(0) +\sum_{n\geq 1  \emph{ odd}} a_n \sqrt{2}
 \sin\left(\frac{\pi nx}{2}\right),
\]
where
\[
a_n = \int_\Gamma f(x) \sqrt{2}  \sin\left(\frac{\pi nx}{2}\right) dx.
\]
\end{thm}

\begin{proof}
For $f \in \textrm{S}([0,1])$, subtract $f(0)$ if necessary so that we may assume
$f(0)=0$. Define
\[
\stackrel{\sim}{f}(x) =
\begin{cases}
f(4x), & 0 \leq x \leq \frac{1}{4} \\
f(2-4x), &  \frac{1}{4} \leq x \leq \frac{1}{2} \\
-f(4x-2), &  \frac{1}{2} \leq x \leq \frac{3}{4} \\
-f(4-4x), &  \frac{3}{4} \leq x \leq 1.
\end{cases}
\]
Then $\stackrel{\sim}{f}$ is piecewise smooth and periodic with period~$1$, so we may
consider it as a function on the circle. Theorem~\ref{StandardFourierTheorem} now
applies, and the Fourier coefficients are easily calculated to be:
\begin{align*}
\tilde{a}_0 &= \int_0^1 \stackrel{\sim}{f}(x) dx =0,\\
\tilde{a}_n &= \int_0^1 \stackrel{\sim}{f}(x)\sqrt{2} \cos(2\pi nx)dx=0, \\
\tilde{b}_n &= \int_0^1 \stackrel{\sim}{f}(x)\sqrt{2} \sin(2\pi nx)dx =
\begin{cases}
 \int_0^1 f(x) \sqrt{2} \sin \left( \frac{\pi nx}{2} \right) dx, & \textrm{if $n$ is odd} \\
\; 0, & \textrm{if $n$ is even.}
\end{cases}
\end{align*}
We can now represent
$\stackrel{\sim}{f}$ by a uniformly convergent series, which in turn gives
a representation of $f(4x)$ for $x \in [0,\frac{1}{4}]$:
\[
f(4x) = \sum_{n\geq 1 \textrm{ odd}} \tilde{b}_n \sqrt{2} \sin(2\pi nx).
\]
Replacing $4x$ with $x$ gives precisely the result we want on $[0,1]$.
\end{proof}

As an application, we prove the following irresistible identity, which was mentioned in
\S\ref{discussion}.

\begin{thm} \label{CoolIdentity}
For all real numbers $0 \leq x,y \leq 1$, we have
\begin{equation*}
\label{eqn:CoolIdentity}
\min\{x,y\} = 8 \sum_{n \geq 1 \emph{ odd}} \frac{\sin \left(\frac{\pi
      nx}{2}\right)
\sin \left(\frac{\pi ny}{2}\right)}{ \pi^2 n^2}.
\end{equation*}
\end{thm}


\begin{proof}
We provide two proofs of this result; the first one is quicker, but the
second proof generalizes better and uses more explicitly the theory of
metrized graphs.

{\bf First proof:}
Fix $y \in [0,1]$ and set $f(x)= \min\{x,y\}$.
Using Theorem~\ref{SpecificFourierTheorem}, we compute that
\begin{align*}
a_n &= \sqrt{2} \int_0^1 f(x) \sin \left(\frac{\pi nx}{2}\right) \, dx \\
&=\sqrt{2}\left\{ \int_0^y x \sin \left(\frac{\pi nx}{2}\right) \, dx + \int_y^1 y\sin
\left(\frac{\pi nx}{2}\right) \, dx \right\} = \frac{4\sqrt{2}}{\pi^2n^2} \sin
\left(\frac{\pi ny}{2}\right).
\end{align*}
Noting that $f(0)=0$ and inserting the coefficients $a_n$ into
Theorem~\ref{SpecificFourierTheorem} yields the result.

{\bf Second proof:}
Thinking of $[0,1]$ as a metrized graph, we see that $\min\{x,y\}$
coincides with the function $j_0(x,y)$ (they have the same Laplacian
and agree at $0$). For $n\geq 1$ odd,
let $\phi_n(x) = \sqrt{2}\sin(\pi nx/2)$ and set $\lambda_n = \pi^2 n^2/4$.

To prove the result, fix $y \in [0,1]$ and use Theorem~\ref{SpecificFourierTheorem} to
write
\begin{equation*} \tag{$\ast \ast$} \label{eq:SpecificFourier2}
j_0(x,y) = \sum_{n\geq 1 \textrm{ odd}} a_n \phi_n(x).
\]

Each $\phi_n$, being an eigenfunction of the Laplacian, satisfies
\[
\Delta \phi_n = \lambda_n \left( \phi_n(x) dx - C_n \delta_0 \right)
\]
for some $C_n \in \RR$ (cf. Definition~\ref{Eigenfunction}). It follows that
\[
\phi_n(x) dx = \frac{\Delta \phi_n}{\lambda_n} + C_n \delta_0.
\]
Applying this to calculate the Fourier coefficients of $j_0(x,y)$, we have
\begin{align*}
a_n &= \int_\Gamma j_0(x,y) \phi_n(x) dx
    = \int_\Gamma j_0(x,y)  \left( \frac{\Delta \phi_n(x)}{\lambda_n} + C_n
    \delta_0(x) \right) \\
    &= \left( \int_\Gamma
    \frac{\phi_n(x)}{\lambda_n} \Delta_x j_0(x,y) \right) + C_n j_0(0,y) \qquad \textrm{by Theorem~\ref{Laplacian Fact 1}}\\
    &= \int_\Gamma \frac{\phi_n(x)}{\lambda_n} \left\{ \delta_y(x) -
    \delta_0(x) \right\} \;
    = \;\frac{\phi_n(y)}{\lambda_n}.
\end{align*}

Substituting our formula for $a_n$ into \eqref{eq:SpecificFourier2}, we obtain
\[
j_0(x,y) = \sum_{n\geq 1  \textrm{ odd}} \frac{\phi_n(x) \phi_n(y)}{\lambda_n},
\]
which is equivalent to the desired result.
\end{proof}

Theorem~\ref{GeneralizedFourierTheorem}, together with an argument
similar to the second proof of Theorem~\ref{CoolIdentity}, yields the following more
general fact:

\begin{thm} \label{IdentityGenerator}
Let $\Gamma$ be a metrized graph, and let $z \in \Gamma$.  Suppose $\{ \phi_n(z) \}_{n
\geq 1}$ is a complete orthonormal set of eigenfunctions of the Laplacian relative to the
measure $\delta_z$, with corresponding eigenvalues $\lambda_n$.  Then for all $x,y \in
\Gamma$, we have
\[
 j_z(x,y) = \sum_{n\geq 1} \frac{\phi_n(x)\phi_n(y)}{\lambda_n}.
\]
\end{thm}

A proof of this theorem and many more results concerning Fourier analysis on metrized graphs can
be found in \cite{BR}.

\begin{exercise}
Find other nice identities like the one in Theorem~\ref{CoolIdentity} by taking a
metrized graph $\Gamma$, working out the eigenfunctions of the Laplacian with respect to
$\delta_z$ for some point $z \in \Gamma$, and applying
Theorem~\ref{IdentityGenerator}.
\end{exercise}

\section{Epilogue}
The material in this expository paper was adapted from a series of lectures
given by the first author and Robert Rumely for the
``Analysis on Metrized Graphs'' REU in summer 2003, and
represents the jumping-off point for several research questions
explored during the REU.  These questions included:
\begin{itemize}
\item Is there a good discrete analogue of the canonical measure? ({\it Yes.})
\item Do the eigenvalues of the Laplacian matrix on a sequence of
  models for a metrized graph $\Gamma$ converge (under suitable hypotheses) to the eigenvalues of
  the Laplacian operator on $\Gamma$? ({\it Yes.})
\item If $\Gamma$ is normalized to have total length 1, can the quantity
$\tau(\Gamma)$ from Exercise~\ref{TauExercise} be arbitrarily small? ({\it No.})
\end{itemize}
A more detailed discussion of these questions, and of the results obtained, can be found
at \verb+http://www.math.uga.edu/~mbaker/REU/REU.html+

\end{document}